\documentclass[12pt]{article}
\usepackage[T2B]{fontenc}
\usepackage{indentfirst}
\usepackage[pdftex, colorlinks, urlcolor=blue, pdfborder={0 0 0 [0 0]}]{hyperref}
\usepackage{amsfonts,amsmath,amssymb,theorem}
\usepackage{array}

\usepackage{esint}
\usepackage{graphicx}
\newtheorem{thm}{Theorem}

\newtheorem{exm}{Example}
\newtheorem{remark}{Remark}
\newtheorem{dfn}{Definition}
\theoremstyle{remark}
\theoremstyle{definition}

\title{The construction of the program control with probability one for stochastic dynamic systems with jumps}

\author{Elena~Karachanskaya
}

\begin{document}
\maketitle

\begin{abstract}
Investigate the stochastic dynamic non-linear system with the Wiener and the Poisson perturbations:
\begin{equation}\label{Puas1-vec2}
\begin{array}{c}
  d {\bf x}(t)= \Bigl( P(t;{\bf x}(t)) +  Q(t;{\bf x}(t)) \cdot {\bf
s}(t;{\bf x}(t))\Bigr) dt + \\
+ \varepsilon\Bigl(t;{\bf x}(t),d {\bf
w}(t),\nu(dt;d\gamma)\Bigr),
\end{array}
\end{equation}
where   ${\bf x}\in \mathbb{R}^{n}$, $n\geq 2$; $\varepsilon(\cdot)$ -- is a strong random perturbation; ${\bf w} (t)$ -- is the $m$-dimensional Wiener process;  $\nu(t;\Delta\gamma)$ --  is the homogeneous on $t$ non-centered Poisson measure.
For such systems we construct the program control ${\bf
s}(t;{\bf x}(t))$ with probability one, which allows the system \eqref{Puas1-vec2} to move on the given trajectory:
$\label{m-fold}
u_{i}(t;{\bf x}(t))=0, \ \ \ i=\overline{1,k}.
$
Control program ${\bf s}(t;{\bf x}(t))$ is  solution of the algebraic system of linear equations. Considered algorithm is based on the first integral theory for stochastic differential equations system.\footnote{%
60H15 (Primary), 37H10 (Secondary)\\
Keywords:
Program control, stochastic dynamic system, Poisson jumps}
\end{abstract}


\section{Preliminary definitions and results}

Suppose that a dynamic system is affected by random disturbance, which takes the form of Wiener and Poisson perturbations. In such  situation we require that important system properties are kept invariant. We investigate the program control construction problem with probability one for the similar systems. We use the first integral concept for the It\^{o}'s stochastic differential equation system (SDES) \cite{D_78,D_89,D_02} and the automorphic function construction \cite{D_04,11_KchUpr}.

Let us consider the random  process $\bf{x}(t)=\bf{x}(t,\omega)$,   $\bf{x}: [0;+\infty) \to  \mathbb{R}^{n}_{+}$ which is the solution of the Cauchy problem for system of the It\^{o}'s generation stochastic differential equations (GSPE)
\begin{equation}\label{1}
\begin{array}{c}
  dx_{i}(t)=a_{i}(t;{\bf{x}}(t))dt+\displaystyle\sum\limits_{k=1}^{m}b_{ik}(t;{\bf{x}}(t))
  dw_{k}(t)
  +
  \displaystyle\int\limits_{R(\gamma)}g_{i}(t;{\bf{x}}(t);\gamma)\nu(dt;d\gamma)\\
{\bf{x}}(t)={\bf{x}}(t,{\bf{x}}_{o})\Bigl|_{t=0}={\bf{x}}_{0},\ \
\ i=\overline{1,n}, \ \ t\geq 0,
\end{array}
\end{equation}
where  ${\bf{w}}(t))$ is the  $m$-dimensional Wiener process, $\nu(t;\Delta\gamma)$ which  is homogeneous on $t$ non centered Poisson measure. This equation we present in the vector form
\begin{equation*}
\begin{array}{c}
  d {\bf x}(t)= A(t;{\bf x}(t)) dt +  B(t;{\bf x}(t))  d {\bf
w}(t)
  +\displaystyle\int\limits_{R(\gamma)}\nu(dt;d\gamma)\cdot
G(t;{\bf x}(t);\gamma)
\end{array}
\end{equation*}

We suppose that the coefficients  ${ a}(t;{\bf x})$, ${ b}(t;{\bf x})$,   and ${g}(t;{\bf x};\gamma)$  are bounded and continuous together with their own derivations  $
\displaystyle \frac{\partial a_{i}(t;{\bf x})}{\partial x_{j}}$,\,
$\displaystyle \frac{\partial b_{i,k}(t;{\bf x})}{\partial
x_{j}}$,\, $\displaystyle \frac{\partial g_{i}(t;{\bf
x};\gamma)}{\partial x_{j}}$,\, $i,j=\overline{1,n}
$   respect to  $(t;{\bf x};\gamma)$, and they are satisfied the conditions of the theorem of existence and uniqueness of the Eq.~(\ref{1}).

In the article \cite{D_78} the concept of the first integral for It\^{o}'s GSDE  system (without Poisson part) has been introduced. In the article \cite{D_02} the concept of the stochastic first integral for It\^{o}'s GSDE  system with the centered measure has been given. Let's introduce the similar concept for the non-centered  Poisson measure.

\begin{dfn}\label{Ayddf1}
Let  random function  $u(t;{\bf x};\omega)$   and   solution ${\bf x}(t)$ of the It\^{o}'s GSDE  system \eqref{1} be defined  on the same    probability space. The  function  $u(t;{\bf x};\omega)$  is called the stochastic first integral for It\^{o}'s GSDE  system \eqref{1} with  non centered Poisson measure, if condition
$$
u\Bigl(t;{\bf x}(t; {\bf x}(0));\omega\Bigr)=u\Bigl(0;{\bf
x}(0);\omega\Bigr)
$$
is held with probability equaled to 1 for the each solution  ${\bf x}(t;{\bf x}(0);\omega)$  of the system \eqref{1}.
\end{dfn}

The random function $u\left(t;{\bf x}(t)\right) \in \mathcal{C}_{t,x}^{1,2}$ defined over the same  probability space as the  random process ${\bf x}(t)$ which is the solution for the system \eqref{1}  is  the first integral of the system \eqref{1} iff the function $u\left(t;{\bf x}(t)\right)$ satisfies the terms  $\left.\mathcal{L}\right)$\label{uslL1} \cite{11_KchItoEng}:
\begin{enumerate}
    \item $b_{i\,k}(t;{\bf x})\displaystyle\frac{\partial u(t;{\bf x})}{\partial
x_{i}}=0$, for all $k=\overline{1,m}$ (compensation of the Wiener's perturbations);
    \item $\displaystyle\frac{\partial u(t;{\bf x})}{\partial
t}+\displaystyle\frac{\partial u(t;{\bf x})}{\partial
x_{i}}\Bigl[a_{i}(t;{\bf x})-
\displaystyle\frac{1}{2}\,b_{j\,k}(t;{\bf x})\frac{\partial
b_{i\,k}(t;{\bf x})}{\partial x_{j}}\Bigr]=0$ (independence of the time);
     \item $u(t;{\bf x})-u\Bigl(t;{\bf x}+g(t;{\bf
    x};\gamma)\Bigr)=0$ for any $\gamma\in R(\gamma)$
         in the whole field of the process definition   (compensation of the Poisson's jumps).
\end{enumerate}

\begin{remark}
If we analyze the concrete realization then the non random function  $u(t;{\bf x})$ is the determinate first integral of the stochastic system.
\end{remark}
The concept of the stochastic first integral for the centered Poisson measure was introduced in \cite{D_02}. The obtained conditions for its realization take the necessity of determining the   intensive Poisson distribution density. In this research that condition is missing. Thus, it makes no difference what is the probability distribution of intensities of Poisson jumps. This case is very important for constructing program controls \cite{11_KchUpr}.
In case a Wiener disturbances only, we construct the program control with probability one in article \cite{08_ChContrEn} on basis the research \cite{D_89}.

\begin{thm}\label{teor-1}{\rm\cite{11_KchUpr}}
Let the determinate function $u(t,{\bf x})$ be continuous together with their own  derivatives with respect to all variables $u(t,{\bf x})$ and the random function  $u(t,{\bf x};\omega)$  defined  over the same  probability space as the  random process ${\bf x}(t)$ which is the solution for system
\begin{equation}\label{1m}
\begin{array}{l}
  d{\bf{x}}(t)=A(t;{\bf{x}}(t))dt+B(t;{\bf{x}}(t))
  d\mathbf{w}(t)  +
  \displaystyle\int\limits_{R(\gamma)}\nu(dt;d\gamma)\cdot
G(t;{\bf x}(t);\gamma)\\
{\bf{x}}(t)={\bf{x}}(t,{\bf{x}}_{o})\Bigl|_{t=0}={\bf{x}}_{0},\ \
\ t\geq 0,
\end{array}
\end{equation}
where ${\bf{x}}\in \mathbb{R}^{n}$, $n\geq 2$;  ${\bf{w}}(t)$ is
an  $m$-dimensional Wiener process, $\nu(t;\Delta\gamma)$ --  is homogeneous on $t$ non centered Poisson measure. Assume the vectors $\vec{e}_{o}$, $\vec{e}_{1}$, ..., $\vec{e}_{n}$ consist an orthogonal basis in $\mathbb{R}_{+}\times\mathbb{R}^{n}$. If function $u(t,{\bf
x};\omega)$ is  the first integral of the system
\eqref{1m}, then the coefficients of an Eq. \eqref{1m} and the function
$u(t,{\bf x})$	are related by conditions:
\begin{enumerate}
    \item coefficients
$ B_{k}(t;{\bf{x}})=\displaystyle\sum\limits_{i=1}^{n}b_{ik}(t;{\bf{x}})\vec{e}_{i}$
    $(k=\overline{1,m})$,  are columns of the matrix
   $B(t;{\bf{x}})$, that belongs to a set of functions
\begin{equation}\label{B}
B_{k}(t;{\bf x})\in \left\{
\begin{array}{c}
q_{oo}(t;{\bf x})\cdot \det \left(
\begin{array}{ccc}
  \vec{e}_{1} & \ldots & \vec{e}_{n} \\
  \displaystyle\frac{\partial u(t; {\bf{x}})}{\partial x_{1}}  & \ldots &
  \displaystyle\frac{\partial u(t; {\bf{x}})}{\partial x_{n}}  \\
  f_{31} & \ldots & f_{3n}\\
  \ldots & \ldots & \ldots  \\
  f_{n1} & \ldots & f_{nn}
\end{array}
\right)
\end{array}
\right\},
\end{equation}
where $q_{oo}(t;{\bf x})$ is an arbitrary non-vanishing function;
    \item coefficient
    $A(t;{\bf{x}})$
   belongs to a set of the functions, defined by
\begin{equation}\label{A}
{{A}}(t;{\bf x})\in \left\{ R(t;{\bf x}) + \displaystyle
\frac{1}{2}\, \sum_{k=1}^{n}\displaystyle \biggl[\frac{\partial
B_{k}(\cdot)}{\partial {\bf x}}\biggr]\cdot B_{k}(\cdot)\right\},
\end{equation}
where
a column matrix $ R(t;{\bf x})$ with components $r_{i}(t;{\bf
x}) $, $i=\overline{1,n}$, are defined by:
$$ C^{-1}(t;{\bf x}) \cdot \det H(t;\mathbf{x})= \vec{e}_{0}+
\displaystyle \sum\limits_{i=1}^{n} r_{i}(t;{\bf x})\vec{e}_{i};
$$
$C(t;{\bf x}) $ is an algebraic adjunct of the element
$\vec{e}_{0}$ of a matrix $H(t;\mathbf{x})$ and $\det C(t;{\bf x})\neq 0$;
a matrix $H(t;\mathbf{x})$ is defined as
\begin{equation}\label{H} H(t;\mathbf{x})=\left[
\begin{array}{cccc}
  \vec{e}_{0} & \vec{e}_{1} & \ldots & \vec{e}_{n} \\
 \displaystyle\frac{\partial u(t;\mathbf{x})}{\partial t} &
 \displaystyle\frac{\partial u(t;\mathbf{x})}{\partial x_{1}} & \ldots &
 \displaystyle\frac{\partial u(t;\mathbf{x})}{\partial x_{n}} \\
h_{30} & h_{31} & \ldots & h_{3n} \\
  \ldots & \ldots & \ldots & \ldots \\
  h_{n+1,0} & h_{n+1,1} & \ldots & h_{n+1,n}
\end{array}
\right],
\end{equation}
and
$\displaystyle \biggl[\frac{\partial B_{k}(t;{\bf x})}{\partial
{\bf x}}\biggr]$ is a matrix of Jacobi for function
$B_{k}(t;{\bf x})$;
    \item coefficient $G(t;{\bf x};\gamma)=
    \displaystyle\sum\limits_{i=1}^{n}g_{i}(t;{\bf{x}};\gamma)\vec{e}_{i}$
     belonging to Poisson measure, is defined by the next representation
     $G(t;{\bf x};\gamma)=\mathbf{y}(t;{\bf x};\gamma)-{\bf x}$,
where $\mathbf{y}(t;{\bf x};\gamma)$ is the solution of the differential equations system
\begin{equation}\label{Y}
\begin{array}{c}
\displaystyle\frac{\partial \mathbf{y}(\cdot;\gamma)}{\partial
\gamma}= \\
=\det\left[
\begin{array}{cccc}
  \vec{e}_{1}& \vec{e}_{2} & \cdots & \vec{e}_{n} \\
  \displaystyle\frac{\partial u(t; \mathbf{y}(\cdot;\gamma))}{\partial y_{1}} &
  \displaystyle\frac{\partial u(t; \mathbf{y}(\cdot;\gamma))}{\partial y_{2}} & \cdots &
  \displaystyle\frac{\partial u(t; \mathbf{y}(\cdot;\gamma))}{\partial y_{n}} \\
  \varphi_{31}(t;\mathbf{y}(\cdot;\gamma)) &
  \varphi_{32}(t;\mathbf{y}(\cdot;\gamma)) & \cdots &
  \varphi_{3n}(t;\mathbf{y}(\cdot;\gamma)) \\
  \cdots & \cdots & \cdots & \cdots \\
   \varphi_{n1}(t;\mathbf{y}(\cdot;\gamma)) &
  \varphi_{n2}(t;\mathbf{y}(\cdot;\gamma)) & \cdots &
  \varphi_{nn}(t;\mathbf{y}(\cdot;\gamma))
\end{array}
\right].
\end{array}
\end{equation}
This solution satisfies the initial conditions:
$\mathbf{y}(t;\mathbf{x};\gamma)\Bigl|_{\gamma=0}=\mathbf{x}$.
\end{enumerate}
The arbitrary functions $f_{ij}=f_{ij}(t,\mathbf{x} )$,
$h_{ij}=h_{ij}(t,\mathbf{x} )$,
$\varphi_{ij}=\varphi_{ij}(t;\mathbf{y}(\cdot;\gamma))$ are defined by the equalities
$
f_{ij}(t,\mathbf{x} )=\displaystyle\frac{\partial
f_{i}(t,\mathbf{x} )}{\partial x_{j}}$, \,
$h_{ij}(t,\mathbf{x} )=\displaystyle\frac{\partial
h_{i}(t,\mathbf{x} )}{\partial x_{j}}$, \,
$\varphi_{ij}(t;\mathbf{y}(\cdot;\gamma))=\displaystyle\frac{\partial
\varphi_{i}(t;\mathbf{y}(\cdot;\gamma))}{\partial y_{j}}
$,
and sets of the functions
$\Bigl\{f_{i}\Bigr\}$, $\Bigl\{h_{i}\Bigr\}$,
$\Bigl\{\varphi_{i}\Bigr\}$ and the function $u(t;\mathbf{x})$ together consist of the class of  independent functions.
\end{thm}

\noindent
{\bf Proof.}
We consider three steps for the proof.

${\bf 1}.$ Let us use the first statement from the conditions
$\left.\mathcal{L}\right)$:
\begin{equation}\label{a}
\displaystyle\sum\limits_{i=1}^{n}b_{i\,k}(t;{\bf
x})\displaystyle\frac{\partial u(t;{\bf x})}{\partial x_{i}}=0, \ \ \ \
\textrm{for all} \ \ k=\overline{1,m}.
\end{equation}
 If
$B_{k}(t;{\bf{x}})=\displaystyle\sum\limits_{i=1}^{n}b_{ik}(t;{\bf{x}})\vec{e}_{i}$
and $\nabla_{\mathbf{x}}u(t;{\bf
x})=\displaystyle\sum\limits_{i=1}^{n}\frac{\partial u(t;{\bf
x})}{\partial x_{i}}\, \vec{e}_{i} $ hold, then Eq.~\eqref{a}  is an orthogonal property of vectors
 $B_{k}(t;{\bf{x}})$ and $\nabla_{\mathbf{x}}u(t;{\bf
x})$.

Taking into account a vector product definition in
$\mathbb{R}^{n}$ and their properties we set a representation for columns
 $B_{k}(t;{\bf{x}})$ of matrix
$B(\cdot)=\Bigl(B_{1}(\cdot),\ldots,B_{m}(\cdot)\Bigr)$:
$$
B_{k}(t;{\bf{x}})\in \left\{
\begin{array}{c}
q_{oo}(t;{\bf{x}})\cdot \det \left(
\begin{array}{ccc}
  \vec{e}_{1} & \ldots & \vec{e}_{n} \\
  \displaystyle\frac{\partial u(t; {\bf{x}})}{\partial x_{1}}  & \ldots &
  \displaystyle\frac{\partial u(t; {\bf{x}})}{\partial x_{n}}  \\
  f_{31} & \ldots & f_{3n}\\
  \ldots & \ldots & \ldots  \\
  f_{n1} & \ldots & f_{nn}
\end{array}
\right)
\end{array}
\right\}; \eqno(\ref{B})
$$
where $f_{i}=f_{i}(t,\mathbf{x} )$, \, $i=\overline{3,n}$,\,
and functions $f_{ij}(t,\mathbf{x} )=\displaystyle\frac{\partial
f_{i}(t,\mathbf{x} )}{\partial x_{j}}$ and
 $u(t,\mathbf{x} )$ together generate the collection of a independence ones.

${\bf 2}.$ Later we use the second statement of the conditions
$\left.\mathcal{L}\right)$:
\begin{equation}\label{b}
\displaystyle\frac{\partial u(t;{\bf x})}{\partial
t}+\displaystyle\sum\limits_{i=1}^{n}\frac{\partial u(t;{\bf
x})}{\partial x_{i}}\Bigl[a_{i}(t;{\bf x})-
\displaystyle\frac{1}{2}\sum\limits_{k=1}^{m}\sum\limits_{j=1}^{n}\,b_{j\,k}(t;{\bf
x})\frac{\partial b_{i\,k}(t;{\bf x})}{\partial x_{j}}\Bigr]=0.
\end{equation}
Let us suppose
$$Q(t;{\bf x})=1+
\displaystyle\sum\limits_{i=1}^{n}a_{i}(t;{\bf{x}})-
\displaystyle\frac{1}{2}\sum\limits_{i=1}^{n}\sum\limits_{k=1}^{m}\sum\limits_{j=1}^{n}\,b_{j\,k}(t;{\bf
x})\frac{\partial b_{i\,k}(t;{\bf x})}{\partial x_{j}}.
$$
According to the scheme in article {\rm{\cite{D_89}}} we consider two vectors:
$$\square
u(t;{\bf x})= \displaystyle\frac{\partial u(t;{\bf x})}{\partial
t}\, \vec{e}_{0}+\sum\limits_{i=1}^{n}\frac{\partial u(t;{\bf
x})}{\partial x_{i}}\, \vec{e}_{i}
$$
and
$$\overrightarrow{Q}(t;{\bf x})=\vec{e}_{0}+
\displaystyle\sum\limits_{i=1}^{n}a_{i}(t;{\bf{x}})\vec{e}_{i}-
\displaystyle\frac{1}{2}\sum\limits_{i=1}^{n}\sum\limits_{k=1}^{m}\sum\limits_{j=1}^{n}\,b_{j\,k}(t;{\bf
x})\frac{\partial b_{i\,k}(t;{\bf x})}{\partial
x_{j}}\,\vec{e}_{i}.
$$

Then the Eq.~\eqref{b} means that the vectors  $\square u(t;{\bf x})$ and
$\overrightarrow{Q}(t;{\bf x})$ are orthogonal. Later we use the vector product and it's properties again, and we obtain:
\begin{equation}\label{Q}
\overrightarrow{Q}(t;{\bf x})\in \left\{\det \left[
\begin{array}{cccc}
  \vec{e}_{0} & \vec{e}_{1} & \ldots & \vec{e}_{n} \\
 \displaystyle\frac{\partial u(t;\mathbf{x})}{\partial t} &
 \displaystyle\frac{\partial u(t;\mathbf{x})}{\partial x_{1}} & \ldots &
 \displaystyle\frac{\partial u(t;\mathbf{x})}{\partial x_{n}} \\
h_{30} & h_{31} & \ldots & h_{3n} \\
  \ldots & \ldots & \ldots & \ldots \\
  h_{n+1,0} & h_{n+1,1} & \ldots & h_{n+1,n}
\end{array}
\right] \right\}= \left\{\det H \right\},
\end{equation}
where  $h_{i}=f_{i}(t,\mathbf{x} )$, \,
$i=\overline{3,n+1}$,\, and functions  $h_{ij}(t,\mathbf{x}
)=\displaystyle\frac{\partial h_{i}(t,\mathbf{x} )}{\partial
x_{j}}$ and $u(t,\mathbf{x} )$  together generate the collection of a independence ones.
Let us consider the vector
$$
\overrightarrow{{\widetilde{A}}}(t;{\bf x})= \vec{e}_{o}+
\displaystyle\sum\limits_{i=1}^{n}a_{i}(t;{\bf x})\vec{e}_{i}=
\overrightarrow{Q}(t;{\bf x}) +
\displaystyle\frac{1}{2}\displaystyle\sum\limits_{i=1}^{n}\sum\limits_{k=1}^{m}\sum\limits_{j=1}^{n}\,b_{j\,k}(t;{\bf
x})\frac{\partial b_{i\,k}(t;{\bf x})}{\partial
x_{j}}\,\vec{e}_{i}.
$$
As a coefficient of  \, $\vec{e}_{o}$ is equal to 1, then we get:
\begin{equation*}
\overrightarrow{{\widetilde{A}}}(t;{\bf x})\in \left\{
C^{-1}(t;{\bf x}) \cdot \det H(t;{\bf{x}}) + \displaystyle \frac{1}{2}\,
\displaystyle\sum\limits_{i=1}^{n}\sum\limits_{k=1}^{m}\sum\limits_{j=1}^{n}\,b_{j\,k}(t;{\bf
x})\frac{\partial b_{i\,k}(t;{\bf x})}{\partial
x_{j}}\,\vec{e}_{i}\right\},
\end{equation*}
where $C(t;{\bf x}) $
is an algebraical adjunct of element
$\vec{e}_{0}$ of matrix $H(t;{\bf{x}})$, $\det C(t;{\bf x})\neq 0 $. As far as
the vector $C^{-1}(t;{\bf x}) \cdot \det H(t;{\bf{x}}) $ we present in a form:
$$ C^{-1}(t;{\bf x}) \cdot \det H(t;{\bf{x}})= \vec{e}_{0}+
\displaystyle \sum\limits_{i=1}^{n} r_{i}(t;{\bf x})\vec{e}_{i},
$$
let us introduce the following vector $ R(t;{\bf x})$ that has components
$r_{i}(t;{\bf x}) $, $i=\overline{1,n}$.

We might the next representation:
\begin{equation}\label{Bk}
\displaystyle\sum\limits_{i=1}^{n}\sum\limits_{k=1}^{m}\sum\limits_{j=1}^{n}\,b_{j\,k}(t;{\bf
x})\frac{\partial b_{i\,k}(t;{\bf x})}{\partial
x_{j}}=\sum_{k=1}^{n}\displaystyle \biggl[\frac{\partial
B_{k}(t;{\bf x})}{\partial {\bf x}}\biggr]\cdot B_{k}(t;{\bf x}),
\end{equation}
where $\displaystyle \biggl[\frac{\partial B_{k}(t;{\bf
x})}{\partial {\bf x}}\biggr]$ is a matrix of Jacobi for function $B_{k}(t;{\bf x})$.
It follows that a coefficient
${{A}}(t;{\bf x})$ is a sum of matrices:
\begin{equation*}
{{A}}(t;{\bf x})\in \left\{ R(t;{\bf x})+ \displaystyle
\frac{1}{2}\, \sum\limits_{k=1}^{m}\sum_{k=1}^{n}\displaystyle
\biggl[\frac{\partial B_{k}(t;{\bf x})}{\partial {\bf
x}}\biggr]\cdot B_{k}(t;{\bf x})\right\},
\end{equation*}

${\bf 3}.$ According to the third statement in
$\left.\mathcal{L}\right)$ for all $\gamma\in R(\gamma)$
we have:
\begin{equation}\label{c}
u(t;{\bf x};\omega)-u\Bigl(t;{\bf x}+G(t;{\bf
x};\gamma);\omega\Bigr)=0.
\end{equation}
It means that the function $u(t;{\bf x};\omega)$ is automorphic function
under translation by
${\bf x}$ with function $G(t;{\bf x};\gamma)$. Let us set a  condition for it.

According to {\rm{\cite{D_02}}} we have ${\bf y}(t;{\bf
x};\gamma)={\bf x}+G(t;{\bf x};\gamma)$. For notational simplicity  we 	
take off a parameter $\omega$. Then we obtain
\begin{equation}\label{u}
u(t;{\bf x})=u\Bigl(t;{\bf y}(t;{\bf x};\gamma)\Bigr)
\end{equation}
for all $\gamma\in R(\gamma)$. Hence, $\displaystyle
\frac{\partial u(t;{\bf x})}{\partial \gamma}=0$ and it holds:
\begin{equation}\label{p}
\displaystyle \frac{\partial u\Bigl(t;{\bf y}(t;{\bf
x};\gamma)\Bigr)}{\partial
\gamma}\equiv\sum\limits_{i=1}^{n}\frac{\partial u\Bigl(t;{\bf
y}(t;{\bf x};\gamma)\Bigr)}{\partial y_{i}}\frac{\partial
y_{i}(t;{\bf x};\gamma)}{\partial \gamma}=0.
\end{equation}
Eq.~{p} denotes that the vectors $\nabla_{{\bf
y}}u\Bigl(t;{\bf
y}(\cdot;\gamma)\Bigr)=\displaystyle\sum\limits_{i=1}^{n}\frac{\partial
u\Bigl(t;{\bf y}(\cdot;\gamma)\Bigr)}{\partial y_{i}}\,
\vec{e}_{i}$ and $\displaystyle\frac{\partial {\bf y}(\cdot;\gamma)
}{\partial \gamma}=\sum\limits_{i=1}^{n}\frac{\partial
y_{i}(\cdot;\gamma)}{\partial \gamma}\, \vec{e}_{i}$ are orthogonal and they
connected by the relation
\begin{equation}\label{Gam}
\displaystyle\frac{\partial {\bf y}(\cdot;\gamma) }{\partial
\gamma}\in \left\{\det\left[
\begin{array}{ccc}
  \vec{e}_{1} & \ldots & \vec{e}_{n} \\
  \displaystyle \frac{\partial
u\Bigl(t;{\bf y}(\cdot;\gamma)\Bigr)}{\partial y_{1}}& \ldots &
\displaystyle \frac{\partial
u\Bigl(t;{\bf y}(\cdot;\gamma)\Bigr)}{\partial y_{n}} \\
  \varphi_{31} & \ldots& \varphi_{3n} \\
\ldots & \ldots& \ldots \\
\varphi_{n1} & \ldots& \varphi_{nn}
\end{array}
\right]\right\}.
\end{equation}
There are $\varphi_{i}(t;{\bf y})$, $i=\overline{3,n}$ such functions that
 $\varphi_{ij} =\displaystyle \frac{\partial \varphi_{i}(t;{\bf
y})}{\partial y_{j}}$ and $u\Bigl(t;{\bf
y}(\cdot;\gamma)\Bigr)$ together generate the collection of a independence functions.

By virtue of the fact that
${\bf y}(t;{\bf x};\gamma)={\bf x}+G(t;{\bf x};\gamma)$  then the Eq.~\eqref{Gam}
is the differential system, which has a function  ${\bf
y}(\cdot;\gamma)$ as unknown. 	Let us expand the determinate \eqref{Gam}
 into the first row. Then we have $\displaystyle\frac{\partial {\bf
y}(\cdot;\gamma) }{\partial \gamma}= \alpha
\displaystyle\sum\limits_{i=1}^{n}S_{i}({\bf y}(\cdot;\gamma))\,
\vec{e}_{i}$, where $\alpha$ is an arbitrary function, which independent
from~${\bf y}$.

Such away we obtain the next differential system:
\begin{equation}\label{Gam1}
\left\{\begin{array}{c}
 \displaystyle\frac{\partial  y_{1}(\cdot;\gamma) }{\partial
\gamma}= \alpha S_{1}({\bf y}(\cdot;\gamma)),\\
\cdots \\
 \displaystyle\frac{\partial  y_{n}(\cdot;\gamma)
}{\partial \gamma}= \alpha S_{n}({\bf y}(\cdot;\gamma)).
\end{array}
\right.
\end{equation}
Suppose that ${\bf y}(t;{\bf x};\gamma;\theta)$ is a solution of
the Eq.~\eqref{Gam1}, where $\theta$ is a constant vector. 	
By virtue of \eqref{u} holds for all
 $t$,\, ${\bf x}$ and for all $\gamma$, then is
\begin{equation}\label{u1}
u(t;{\bf x})=u\Bigl(t;{\bf y}(t;{\bf
x};\gamma_{1};\theta)\Bigr)=u\Bigl(t;{\bf x}+ G(t;{\bf
x};\gamma_{1};\theta)\Bigr)=u\Bigl(t;{\bf x}+ G(t;{\bf
x};\gamma_{2};\theta)\Bigr).
\end{equation}
As a special case for same  $\gamma=\gamma_{o}$ the Eq.~\eqref{u1}
is determined by expression $G(t;{\bf
x};\gamma_{o};\theta)=0$.

Unconstrained in common we suppose that
$$
G(t;{\bf x};\gamma_{o};\theta)\equiv G(t;{\bf x};0)=0.
$$
Such away, the system \eqref{u} which has an initial conditions
$\mathbf{y}(t;\mathbf{x};\gamma)\Bigl|_{\gamma=0}=\mathbf{x}$ and it has an uniqueness solving. This system we transform to the form:
$$
\displaystyle\frac{\partial \mathbf{y}(\cdot;\gamma)}{\partial
\gamma}= \det\left[
\begin{array}{cccc}
  \vec{e}_{1}& \vec{e}_{2} & \cdots & \vec{e}_{n} \\
  \displaystyle\frac{\partial u(t; \mathbf{y}(\cdot;\gamma))}{\partial y_{1}} &
  \displaystyle\frac{\partial u(t; \mathbf{y}(\cdot;\gamma))}{\partial y_{2}} & \cdots &
  \displaystyle\frac{\partial u(t; \mathbf{y}(\cdot;\gamma))}{\partial y_{n}} \\
  \varphi_{31}(t;\mathbf{y}(\cdot;\gamma)) &
  \varphi_{32}(t;\mathbf{y}(\cdot;\gamma)) & \cdots &
  \varphi_{3n}(t;\mathbf{y}(\cdot;\gamma)) \\
  \cdots & \cdots & \cdots & \cdots \\
   \varphi_{n1}(t;\mathbf{y}(\cdot;\gamma)) &
  \varphi_{n2}(t;\mathbf{y}(\cdot;\gamma)) & \cdots &
  \varphi_{nn}(t;\mathbf{y}(\cdot;\gamma))
\end{array}.
\right] \eqno(\ref{Y})
$$
Consequently, the function $G(t;{\bf x};\gamma)$ is automorphic transform for the function
 $u(t;{\bf x})$ and it has the next representation:
     $G(t;{\bf x};\gamma)=\mathbf{y}(t;{\bf x};\gamma)-{\bf x}$,
where $\mathbf{y}(t;{\bf x};\gamma)$ is a solution of the differential system
\eqref{Y} with an initial condition
$\mathbf{y}(t;\mathbf{x};\gamma)\Bigl|_{\gamma=0}=\mathbf{x}$.

\section{Construction of the program control}

On the analogy of the article \cite{08_ChContrEn} 	let us introduce the following definition.
\begin{dfn}\label{df2}
By the program motion of the stochastic system
\begin{equation}\label{2}
\begin{array}{c}
  d {\bf x}(t)= \Bigl[ P(t;{\bf x}(t)) +  Q(t;{\bf x}(t)) \cdot {\bf
s}(t;{\bf x}(t))  \Bigr] dt + \\
+B(t;{\bf x}(t))  d {\bf
w}(t)
 +\displaystyle\int\limits_{R(\gamma)}\nu(dt;d\gamma)G(t;{\bf
x}(t);\gamma),
\end{array}
\end{equation}
where ${\bf w}(t)$ is the $m$-dimensional Wiener process;
$\nu(t;\triangle \gamma)$ --  is homogeneous on $t$ non centered Poisson measure,
is meant the solution ${\bf x}(t; {\bf x}_{o},{\bf s};\omega)$ which enables one under some (program) control ${\bf s}(t;{\bf x})$ and for all $t $ to remain with the probability 1 on the given integral manifold $
u\Bigl(t;{\bf x}(t;{\bf x}_{o})\Bigr)= u(0;{\bf x}_{o}), $
which is the first integral of the Eq.~\eqref{2} under the given initial conditions
$${\bf x}(t;{\bf x}_{o})\Bigr|_{t=0}={\bf x}_{o}.
$$
\end{dfn}

By this means we can construct the program control with the probability equaled to 1 for dynamic systems which are subjected to the Wiener perturbations and the Poisson jumps.

\begin{thm}\label{teor-2}
The control ${\bf s}(t;\mathbf{x})$ allowing system  \eqref{2} which is subjected to the Wiener perturbations and the Poisson jumps
always remain  on the dynamically structured integral manifold $ u\Bigl(t;{\bf
x}(t;{\bf x}_{o});\omega\Bigr)= u(0;{\bf x}_{o})$ is determined as the solution of the system of linear equation, which consists of the Eq.~\eqref{2} and  the Eq.~\eqref{1m}. The coefficients of the second equation and respectively coefficient of the first equation of this system are determined by the Theorem~\ref{teor-1}. In addition we completely define the response to random action.
\end{thm}

\begin{exm} Needed is to determine the control ${\bf s}(t;\mathbf{x})$ and the response to random action enabling to dynamical system
\begin{equation}
\begin{array}{c}
  \displaystyle\frac{dx_{1}(t)}{dt}= \Bigl(x_{1}(t)+x_{2}(t)+ e^{-t} +s_{1}(t;\mathbf{x}(t))\Bigr)dt+\\
  + b_{1}(t;\mathbf{x}(t))dw(t)+
  \displaystyle\int_{R(\gamma)}{g}_{1}(t;\mathbf{x}(t);\gamma)\nu(dt;d\gamma), \\
   \displaystyle\frac{dx_{1}(t)}{dt}= \Bigl(x_{1}(t)x_{2}(t)+ e^{-2t}+s_{2}(t;\mathbf{x}(t))\Bigr)dt+\\
   + b_{2}(t;\mathbf{x}(t))dw(t)+
  \displaystyle\int_{R(\gamma)}{g}_{2}(t;\mathbf{x}(t);\gamma)\nu(dt;d\gamma),
\end{array}
\end{equation}
to remain on the given integral manifold
 $u(t;\mathbf{x}(t))=x_{2}(t)e^{-2x_{1}(t)}$ with probability one.
\end{exm}

\textbf{Solution.} According to the Theorem~\ref{teor-2} we construct the GSDЕ system, which has a function $u(t;\mathbf{x})=x_{2}e^{-2x_{1}}$ as the first integral.

We first find an automorphic transformation of function  $u(t;\mathbf{x})=x_{2}e^{-2x_{1}}$. The corresponding partial derivatives of function $u(t;\mathbf{y})=y_{2}e^{-2y_{1}}$ are given by
\begin{equation}
\begin{array}{cc}
  \displaystyle\frac{\partial u(\mathbf{y};t)}{\partial
  y_{1}}=-2y_{2}e^{-2y_{1}},
  & \ \ \ \
  \displaystyle\frac{\partial u(\mathbf{y};t)}{\partial
  y_{2}}=e^{-2y_{1}}.
\end{array}
\end{equation}
Then we have
\begin{equation*}
\displaystyle\frac{\partial
\mathbf{y}(t;\mathbf{x};\gamma)}{\partial \gamma}\equiv
\left(\begin{array}{c}
  \displaystyle\frac{\partial y_{1}(t;\mathbf{x};\gamma)}{\partial
\gamma} \\
  \displaystyle\frac{\partial y_{2}(t;\mathbf{x};\gamma)}{\partial
\gamma}\end{array}
 \right)= \left(\begin{array}{c}
         e^{-2y_{1}} \\
          2y_{2}e^{-2y_{1}}
        \end{array}
\right).
\end{equation*}
Solution of this system is:
\begin{equation*}
\begin{array}{c}
  y_{1}(t;\mathbf{x};\gamma)=\displaystyle\frac{1}{2}\ln\left(2\gamma +2C_{1}(\mathbf{x})\right)
\\
y_{2}(t;\mathbf{x};\gamma)=C_{2}(\mathbf{x})\left(\gamma+C_{1}(\mathbf{x})\right).
\end{array}
\end{equation*}
With regard to initial conditions (according to Theorem~\ref{teor-1}, point 3), i. e.
$\mathbf{y}(t;\mathbf{x};\gamma)\Bigl|_{\gamma=0}=\mathbf{x},$
 we obtain:
\begin{equation*}
\begin{array}{cc}
  C_{1}(\mathbf{x})=\displaystyle\frac{1}{2}e^{2x_{1}},  & \ \ \
  C_{2}(\mathbf{x})=2x_{2}e^{-2x_{1}}.
\end{array}
\end{equation*}
And thus we have:
\begin{equation*}
\begin{array}{c}
  y_{1}(t; \mathbf{x};\gamma)=\displaystyle\frac{1}{2}\ln\left(2\gamma +e^{2x_{1}}\right),\\
  y_{2}(t; \mathbf{x};\gamma)=2x_{2}\gamma e^{-2x_{1}}+x_{2}.
\end{array}
\end{equation*}
Consequently, automorphic transformation $g(\cdot)=(g_{1}(\cdot),g_{2}(\cdot))^{*}$ of the function  $u(t;\mathbf{x})=x_{2}e^{-2x_{1}}$
is:
\begin{equation}\label{g1}
\begin{array}{c}
  g_{1}(t; \mathbf{x};\gamma)=\displaystyle\frac{1}{2}\ln\left(2\gamma +e^{2x_{1}}\right)-x_{1}, \\
  g_{2}(t; \mathbf{x};\gamma)=2x_{2}\gamma e^{-2x_{1}}.
\end{array}
\end{equation}
According  to  the statement 1 of the Theorem~\ref{teor-1} we construct a matrix
 $B$ (in this case it is a column as far as
${\bf w}(t)$ is  one-dimensional Wiener process):
$$
B(\cdot)\equiv B (t; \mathbf{x})=q_{oo}\left(e^{-2x_{1}}, 2x_{2}e^{-2x_{1}}\right)^{*},
$$
where $q_{oo}=q_{oo}(t; \mathbf{x})$. Let us define the second expression in Eq.~(\ref{A}):
$$
\biggl[\frac{\partial B(t;{\bf x})}{\partial {\bf
x}}\biggr]=q_{oo}\left(
\begin{array}{cc}
  -2e^{-2x_{1}} & 0\\
  4x_{2}e^{-2x_{1}} & 2e^{-2x_{1}}
\end{array}
\right),
$$
$$
\biggl[\frac{\partial B(t;{\bf x})}{\partial {\bf
x}}\biggr]B(t;{\bf x})=q_{oo}^{2}\left(
\begin{array}{c}
 -4e^{-2x_{1}} \\
0
\end{array}
\right)=\left(
\begin{array}{c}
-4q_{oo}^{2}e^{-4x_{1}} \\
0
\end{array}
\right).
$$

By Eq.~(\ref{H})  let's construct the matrix $H(t;{\bf x})$ and their determinant
$$
\det H(t;{\bf x})=\det\left(
\begin{array}{ccc}
  \vec{e}_{0} & \vec{e}_{1} & \vec{e}_{2} \\
  0 & -2x_{2}e^{-2x_{1}} & e^{-2x_{1}} \\
  f_{1} &  f_{2} &  f_{3}
\end{array}
\right)= $$
$$=
\vec{e}_{0}\Bigl(-2f_{3}x_{2} e^{-2x_{1}}-f_{2} e^{-2x_{1}}\Bigr)+
\vec{e}_{1}\Bigl(f_{1}
e^{-2x_{1}}\Bigr)+$$
$$+\vec{e}_{2}\Bigl(2f_{1}x_{2} e^{-2x_{1}}\Bigr),
$$
where $f_{i}=f_{i}(t; \mathbf{x})$, $i=1,2,3$.

In result the components of vector  $A=A(t; \mathbf{x})$ are given by
$$
\begin{array}{c}
  a_{1}= -\displaystyle\frac{f_{1}}{f_{2}+2f_{3}x_{2}}+2 q_{oo}^{2}e^{-4x_{1}},\\
  a_{2}= -\displaystyle\frac{2f_{1}x_{2}}{f_{2}+2f_{3}x_{2}}.
\end{array}
$$
Hence, the desired GSDE system is:
$$
\begin{array}{c}
  dx_{1}(t)=\left[-\displaystyle\frac{f_{1}}{f_{2}+2f_{3}x_{2}(t)}+2 q_{oo}^{2}e^{-4x_{1}(t)}\right]dt + \\
   +q_{oo}e^{-2x_{1}(t)}dw(t)
  + \displaystyle\int_{R(\gamma)}\Bigl(\displaystyle\frac{1}{2}\ln\left(2\gamma +e^{2x_{1}(t)}\right)-x_{1}(t)\Bigr)\nu(dt;d\gamma)\\
  dx_{2}(t)=\left[-\displaystyle\frac{2f_{1}x_{2}(t)}{f_{2}+2f_{3}x_{2}(t)}\right]dt
  +\\
  +  q_{00}2x_{2}(t)e^{-2x_{1}(t)}dw(t)
  + \displaystyle\int_{R(\gamma)}(2x_{2}(t)\gamma e^{-2x_{1}(t)})  \nu(dt;d\gamma)
\end{array}
$$
Then we must solve the system of linear equation
$$
\begin{array}{c}
  x_{1}(t)+x_{2}(t)+ e^{-t} +s_{1}(t;\mathbf{x}(t))
  =-\displaystyle\frac{f_{1}}{f_{2}+2f_{3}x_{2}(t)}
  +2 q_{oo}^{2}e^{-4x_{1}(t)}, \\
  x_{1}(t)x_{2}(t)+
  e^{-2t}+s_{2}(t;\mathbf{x}(t))
  =-\displaystyle\frac{2f_{1}x_{2}(t)}{f_{2}+2f_{3}x_{2}(t)}.
\end{array}
$$
Hence, program control is of the form of
$$
\begin{array}{c}
  s_{1}(t;\mathbf{x}(t))=-\displaystyle\frac{f_{1}}{f_{2}+2f_{3}x_{2}(t)}+2 q_{oo}^{2}e^{-4x_{1}}-x_{1}(t)
  -x_{2}(t)- e^{-t},  \\
  s_{2}(t;\mathbf{x}(t))=-\displaystyle\frac{2f_{1}x_{2}(t)}{f_{2}+2f_{3}x_{2}(t)}-x_{1}(t)x_{2}(t)-e^{-2t},
\end{array}
$$
where $f_{i}=f_{i}(t; \mathbf{x}(t))$, $i=1,2,3$, $q_{oo}=q_{oo}(t;
\mathbf{x}(t))$ and $f_{2}+2f_{3}x_{2}(t)\neq 0$. 	The response to the  Wiener action is defined as
$$
\begin{array}{c}
  b_{1}(t;\mathbf{x}(t))=q_{oo}e^{-2x_{1}(t)}, \ \ \ \
    b_{2}(t;\mathbf{x}(t))=q_{oo}2x_{2}(t)e^{-2x_{1}(t)}.
\end{array}
$$
Components of the response to the  Poisson jumps are defined as:
$$\begin{array}{c}
  g_{1}(t; \mathbf{x}(t);\gamma)=\displaystyle\frac{1}{2}\ln\left(2\gamma +e^{2x_{1}(t)}\right)-x_{1}(t),\\
  g_{2}(t; \mathbf{x}(t);\gamma)=2x_{2}(t)\gamma e^{-2x_{1}(t)}.
\end{array}
$$

Therefore,  the set of vectors of the program controls and their corresponding responses to a random action enabling the dynamical system under consideration to remain on the given time-variable manifold are determined.

The choice of  functions
 $f_{i}(t; \mathbf{x}(t))$,
$i=1,2,3$ and $q_{oo}(t; \mathbf{x}(t))$ allows to construct the program control on account of some conditions. For example we take into account a simulation utility.

\begin{remark} 	It is worthy to note that the manifold can be defined by a set of functions
{\rm\cite{11_KchItoEng}}.
\end{remark}

\section*{Conclusion}
Application of the first integral stochastic theory allows to construct the program control for dynamic system with probability equaled to 1 of there are strong random perturbations, caused by Wiener and Poisson processes.

\end{document}